\documentclass{amsart}

\usepackage{stmaryrd}

\usepackage{amsmath, amsthm, amssymb}

\theoremstyle{definition}

\usepackage{scalerel}
\usepackage{mathrsfs}

\usepackage{tikz}
\usepackage{hyperref}
\usepackage{placeins}

\usepackage{proof}

\usepackage{verbatim}
\usepackage{mathtools}
\usepackage{cleveref}

\usepackage{bbm}

\usepackage[bottom]{footmisc}

\usepackage{enumitem}

\usepackage{biblatex} 
\newtheorem{theorem}{Theorem}
\newtheorem*{theorem*}{Theorem}
\newtheorem*{maintheorem*}{Main Theorem}

\newtheorem{prop}[theorem]{Proposition}
\newtheorem{lemma}[theorem]{Lemma}

\newtheorem{claim}{Claim}[theorem]

\newenvironment{claimproof}[1]{\par\noindent\textit{Proof of the Claim.}\space#1}{\hfill $\blacksquare$}

\newcommand\A{\mathscr{A}}
\newcommand\I{\mathscr{I}}
\newcommand\G{\mathcal{G}}

\newcommand\ZFUR{\textup{ZFU}_\text{R}}
\newcommand\ZFCUR{\textup{ZFCU}_\text{R}}

\newcommand{\barx}{\bar{x}}
\newcommand{\bary}{\bar{y}}

\def\<#1>{\left\langle#1\right\rangle}
\def\[#1]{\left\llbracket#1\right\rrbracket}
\renewcommand{\restriction}{\mathord{\upharpoonright}}

\title{\textbf{Abstraction Principles and the Size of Reality}}
\author{Bokai Yao}\address[Bokai Yao]{Peking University}
 \email{bkyao@pku.edu.cn}
\urladdr{https://bokaiyao.com}

\usepackage{graphicx}

\begin{document}
\maketitle
\begin{abstract}
The Fregean ontology can be naturally interpreted within set theory with urelements, where objects correspond to sets and urelements, and concepts to classes. Consequently, Fregean abstraction principles can be formulated as set-theoretic principles. We investigate how the size of reality—i.e., the number of urelements—interacts with these principles. We show that Basic Law V implies that for some well-ordered cardinal $\kappa$, there is no set of urelements of size $\kappa$. Building on recent work by Hamkins \cite{hamkins2022fregean}, we show that, under certain additional axioms, Basic Law V holds if and only if the urelements form a set. We construct models of urelement set theory in which the Reflection Principle holds while Hume's Principle fails for sets. Additionally, assuming the consistency of an inaccessible cardinal, we produce a model of Kelley-Morse class theory with urelements that has a global well-ordering but lacks a definable map satisfying Hume's Principle for classes.
\end{abstract}

\section{Introduction}
Fregean abstractions principles aim to establish a certain correspondence between objects and concepts. Frege's \textit{Basic Law V} asserts that each concept is associated with an object, called its \textit{extension}, such that two concepts are co-extensional if and only if they have the same extension. Another prominent instance of Fregean abstraction is \textit{Hume's Principle}, according to which every concept is associated with an object, called its \textit{number}, such that two concepts are equinumerous if and only if they have the same number.  Although Frege's system is inconsistent due to Russell's paradox, the neo-logicist project has succeeded in using consistent fragments of Frege's system to recover arithmetic. On the other hand, the Zermelo-Fraenkel set theory, a different framework from Frege's, provides a unified foundation for mathematics. It is thus natural to investigate the relationship between these two foundational frameworks.

The relationship between abstraction principles and set theory has been studied extensively in the literature. Most of these studies focus on two issues. One is the examination of the model-theoretic properties of theories of abstraction principles using ZFC as a meta-theory (e.g., see Fine \cite{Fine2002-FINTLO-4}, Shapiro \cite{Shapiro2005-SHASAA-12}, and Shapiro $\&$ Roberts \cite{Roberts2023-ROBHPB}). The other is the attempt to develop set theory itself through abstraction principles (e.g., see Boolos \cite{Boolos1989-BOOIA}, Shapiro $\&$ Weir \cite{Shapiro1999-RICNVZ}, and Shapiro $\&$ Uzquiano \cite{Shapiro2008-SHAFMZ-2}). However, there is another perspective to consider: how do abstraction principles behave if we formulate them as set-theoretic principles?

To begin with, the Fregean framework can be naturally interpreted in the iterative conception of set: at the initial stage we have some basic objects called \textit{urelements}, i.e., \textit{non-sets} that are members form sets, and then there are sets of these urelements, sets of sets of them, and so on. In this picture, Fregean concepts can be seen as \textit{classes} of the objects formed in the iterative procedure, and the concepts that are too big to coincide with sets are the \textit{proper classes}. The connection between the Fregean framework and set theory extends beyond this analogy. Recent work by Hamkins \cite{hamkins2022fregean} shows that every model of ZF equipped with definable classes has a second-order definable map fulfilling Basic Law V (\cite[Theorem 1]{hamkins2022fregean}) as well as a second-order definable map fulfilling Hume's Principle (\cite[Theorem 7]{hamkins2022fregean}).

However, ZF set theory is not the most natural set-theoretic framework for studying abstraction principles: every object in ZF is assumed to be a set.  Yet, abstraction principles should be \textit{universally applicable} (Shapiro \cite{Shapiro2005-SHASAA-12} and Shapiro $\&$ Roberts \cite{Roberts2023-ROBHPB}), meaning they are intended to talk about \textit{all} concepts along with \textit{all} objects. Therefore, it is more natural to consider abstraction principles within \textit{set theory with urelements}, a framework that allows all kinds of objects to be members of sets. A natural question immediately arises: How do Fregean abstraction principles behave in such set theory?

We show that abstraction principles in urelement set theory turn out to be contingent upon the size of reality, i.e., how many urelements there are. For example, assuming certain additional axioms, Basic Law V holds if and only if the urelements form a set (Theorem \ref{thm:BLVequivalence1} and \ref{thm:BLVequivalence2} ); and Basic Law V by itself implies that there cannot be unboundedly many urelements (Theorem \ref{plenitude->notBLV}). In the case of Hume's Principle, we improve a result by Gauntt\cite{gauntt1967undefinability} and L\'evy \cite{levy1969definability} by constructing models of ZF in which the urelements forma a proper class and Hume's Principle fails for sets (Theorem \ref{thm:HPindependence}). Furthermore, we construct a model of Kelley-Morse class theory in which the urelements are more numerous than the pure sets. As a result, the model, despite having a global well-ordering, has no definable map fulfilling Hume's Principle for classes (Theorem \ref{KMCUnvdashHP}).

\section{Preliminaries on Urelement Set Theory}
\subsection{Axioms in urelement set theory}\label{section:AxiomsinZFU}
The language of urelement set theory, in addition to $\in$, contains a unary predicate $\A$ for urelements. The axioms of $\ZFUR$ (R for Replacement) include Foundation, Pairing, Union, Powerset, Infinity, Separation, Replacement, Extensionality \textit{for sets}, and the axiom that no urelements have members (see \cite[Section 1]{YAO_2024} for the precise formulation of these axioms). $\ZFCUR$ is $\ZFUR$ + the Axiom of Choice. Note that $\ZFCUR$ allows a proper class of urelements.

Although $\ZFCUR$ looks very much like \textit{the} urelement analog of ZFC, there turns out to be a hierarchy of axioms that are independent of $\ZFCUR$ ( \cite[Theorem 2.5]{YAO_2024}), some of which are ZF-theorems. In particular, the following two set-theoretic principles are not provable from $\ZFCUR$.
\begin{itemize}
\item[] (Collection) $\forall w, u (\forall x \in w \ \exists y \varphi(x, y, u)   \rightarrow \exists v \forall x \in w\  \exists y \in v\  \varphi(x, y, u))$.
\item [] (RP) For every set $x$ there is a transitive set $t$ with $x \subseteq t$ such that for every $x_1, ..., x_n \in t$, $\varphi(x_1, ..., x_n) \leftrightarrow \varphi^t(x_1, ..., x_n)$. 
\end{itemize}
In particular, $\ZFCUR$ has models where the urelements form a proper class but every set of urelements is finite (\cite[Theorem 2.17]{YAO_2024})---call these models \textit{finite-kernel models}. Collection fails in \textit{finite-kernel} models, because there is a set of urelements of size $n$ for every $n \in \omega$ but there is no corresponding collection set. Consequently, RP also fails in finite-kernel models as RP is equivalent to Collection over $\ZFCUR$ (\cite[Theorem 2.5]{YAO_2024}).\footnote{While RP implies Collection without AC, it is not known if Collection implies RP over $\ZFUR$.}

Among the principles that are independent of $\ZFCUR$, Collection and RP are of particular interest because they are often regarded as fundamental axioms of set theory. For example, it is a standard result that Collection is equivalent to Replacement over the remaining axioms of ZF. Consequently, many axiomatizations of ZF (e.g., \cite{chang1990model} and \cite{bell1977boolean}) adopt Collection instead of Replacement as an axiom—hence our use of the subscript R to emphasize that ZF(C)U$_\text{R}$ includes only Replacement. Moreover, certain natural fragments of ZFC, such as ZFC without Powerset, behave nicely only when Collection holds (\cite{Gitman2016-GITWIT}).

RP is a classic version of set-theoretic reflection principles, which articulates the philosophical idea that the set-theoretic universe is \textit{indescribable}—that is, whatever holds true in the universe is already reflected by an initial fragment. This \textit{indescribability conception} of sets is particularly powerful, as it provides a unified justification for many axioms, including Collection, Pairing, and Infinity, all of which follow from RP (\cite{levy1961principles}). In the second-order context, RP is known as Bernays' Reflection Principle, and it is able to justify various large cardinal axioms (see \cite{bernays1976problem} and \cite{tait2005constructing}).

Accordingly, ZF with urelements allows for a certain degree of \textit{axiomatization freedom}, meaning that strong principles such as Collection or RP can be incorporated into its axiomatization. As a result, independent results concerning ZF with urelements become more robust when established within the framework of $\ZFUR$ + RP (or $\ZFUR$ + Collection). We will return to this point in Section \ref{section:HPforSets}.

\subsection{Basic facts in $\ZFUR$}\label{section:BasicFacts}
Now we review some basic notation and facts about $\ZFUR$. Every object $x$ in $\ZFUR$ has a \textit{kernel}, denoted by $ker(x)$, which is the set of the urelements in the transitive closure of $\{x\}$. A set is pure if its kernel is empty. $V$ denotes the class of all pure sets. $Ord$ is the class of all ordinals, which are transitive \textit{pure} sets well-ordered by the membership relation. $\kappa$ is a \textit{well-ordered cardinal} if it is an initial ordinal. $\A$ will also stand for the class of all urelements. $A \subseteq \A$ thus means ``$A$ is a set of urelements''. For any $A \subseteq \A$, the $V_{\alpha}(A)$- hierarchy is defined as usual, i.e.,
\begin{itemize}
    \item [] $V_0(A) = A$;
    \item [] $V_{\alpha+1}(A) = P(V_{\alpha}(A)) \cup A$;
    \item [] $V_{\gamma}(A) = \bigcup_{\alpha < \gamma} V_\alpha(A)$, where $\gamma$ is a limit;
    \item [] $V(A) = \bigcup_{\alpha \in Ord} V_\alpha(A)$.
\end{itemize}
For every $x$ and set $A \subseteq \A$, $x \in V(A)$ if and only if $ker(x) \subseteq A$. $U$ denotes the class of all objects, i.e., $U = \bigcup_{A\subseteq\A} V(A)$. The rank of an object $x$, denoted by $\rho (x)$, is the least ordinal $\alpha$ such that $x \in V_\alpha (A)$ for some $A\subseteq \A$. When there is only a set of urelements, for every $\alpha$ the objects of rank $\alpha$ form a set. An important feature of $\ZFUR$ is that $U$ has many non-trivial automorphisms. Every permutation $\pi$ of a set of urelements can be extended to a definable permutation of $\A$ by letting $\pi$ be identity elsewhere, and $\pi$ can be further extended to a permutation of $U$ by letting $\pi x$ be $\{\pi y : y \in x \}$ for every set $x$. Such $\pi$ preserves $\in$  and thus is an automorphism of $U$. For every $x$ and automorphism $\pi$, $\pi$ point-wise fixes $x$ whenever $\pi$ point-wise fixes $ker(x)$.

We will also discuss \textit{class theories} with urelements. For instance, every model $U$ of $\ZFUR$ equipped with definable classes will be a model of \textit{G\"odel-Benarys class theory} with urelements (GBU), while \textit{Kelley-Morse class theory} with urelements (KMU) is the stronger theory which includes the impredicative comprehension scheme. Every permutation $\pi$ of a set of urelements can also be extended to a permutation of classes that preserves second-order assertions by letting $\pi X = \{\pi x : x \in X\}$ for every class $X$.

\section{Basic law V in Urelement Set Theory}
A map $X \mapsto \epsilon X$ from classes to first-order objects \textit{fulfills Basic Law V} if for any classes $X$ and $Y$, $\epsilon X = \epsilon Y \leftrightarrow \forall x (x \in X \leftrightarrow x \in Y)$. Hamkins \cite{hamkins2022fregean} shows that in every model of ZF equipped with definable classes, there is a second-order definable map that fulfills Basic Law V. We shall first present Hamkins' theorem in $\ZFUR$. The proof uses \textit{Scott's trick}, namely, every class can be represented by the \textit{set} of its elements with minimal rank, which is available as long as the urelements form a set. 

\begin{theorem}[{\cite[Thereom 1]{hamkins2022fregean}}]\label{thm:ASet->BLV}
Let $U$ be a model of $\ZFUR$ + ``$\A$ is a set'' equipped with definable classes. Then $U$ has a definable map that fulfills Basic Law V.
\end{theorem}
\begin{proof}
We first in the meta-theory fix a enumeration $\psi_0, ..., \psi_n, ... $ of the formulas of urelement set theory, and this enumeration will be the standard part of a definable enumeration in $U$. Furthermore, for every standard natural number $k$, $\ZFUR$ has a definable $\Sigma_k$-truth predicate. Thus, given a definable class $X$ of $U$, we let $\varphi(X, \epsilon X)$ be the second-order assertion
\begin{enumerate}
\item [] ``$\epsilon X $ is an ordered pair $\<\ulcorner \psi_n \urcorner,  u>$, where $\psi_n$ is a $\Sigma_k$ formula for some $k$ and $u$ is the set of parameters $p$ with minimal rank such that there is a $\Sigma_k$ truth predicate $T$ with $\forall x (x \in X \leftrightarrow T(\ulcorner \psi_n \urcorner, \<x, p>))$, and no preceding formula $\psi_i$ has this property.''
\end{enumerate}
The map $X \mapsto \epsilon X$ then fulfills Basic Law V (see \cite[page 6]{hamkins2022fregean}).
\end{proof}
\noindent The definability of the map $X \mapsto \epsilon X$, as Hamkins notes, is a major difference between his result and the earlier consistency proofs of Basic Law V  given by Parsons \cite{Parsons1987-TEROTC}, Bell \cite{Bell1994-BELFEO} and Burgess \cite{Burgess1998-JOHOAC}. However, in urelement set theory such definable maps do not always exist.
\begin{lemma}\label{BLV<->Aset}
Let $U$ be a model of $\ZFUR$ in which the following holds.
\begin{enumerate}
\item [] (*) For every $A \subseteq \A$, there is a countably infinite $B \subseteq \A$ such that $B \cap A = \emptyset$.
\end{enumerate}
Then $U$, equipped with definable classes, has no parametrically definable map that fulfills Basic Law V.
\end{lemma}
\begin{proof}
Suppose \textit{for reductio} that in $U$ some second-order formula $\varphi(X, \epsilon X, P)$ with a parameter $P$ defines a extension-assignment map such that for every definable classes $X, Y$ of $U$, 
$$\epsilon X = \epsilon Y \leftrightarrow  \forall x (x 
\in X \leftrightarrow x \in Y).$$
$P$ is a definable class of $U$ so $P = \{ x \in U : U \models \psi (x, z)\}$ for some first-order formula $\psi$ and $z \in U$. In $U$, let $A \subseteq \A$ be an infinite set of urelements disjoint from $ker(z)$ and define $X = \{B \subseteq \A : B - A \text{ is finite} \}$.

\begin{claim}\label{claim1.1}
$A \cap ker(\epsilon X)$ is not empty.
\end{claim}
\begin{claimproof}
Suppose otherwise. Let $A' \subseteq \A$ be such that $A' \sim \omega$ and $A' \cap (A \cup ker(\epsilon X) \cup ker(z)) = \emptyset$, which exists by (*). Let $A_1$ be a countably infinite subset of $A$ and $\pi$ be an automorphism that swaps $A_1$ and $A'$ while point-wise fixing everything else (in particular, $\epsilon X$ and $P$). Since $\varphi (X, \epsilon X, P)$, $\varphi (\pi X, \pi \epsilon X, \pi P)$ and so $\varphi (\pi X, \epsilon X, P)$. Thus, $\pi X$ and $X$ must be co-extensional. But $\pi X = \{ B \subseteq \A: B - \pi A \text{ is finite}\}$, and $A_1 \in X$ but $A_1 \notin \pi X$--contradiction.\end{claimproof}

\vspace{4pt}
\noindent Fix some urelements $a \in ker(\epsilon X) \cap A$ and $b \notin ker(\epsilon X)\cup ker(z)$. Let $\sigma$ be the automorphism that only swaps $a$ and $b$. Consequently, $\sigma \epsilon X \neq \epsilon X$. As $X$ and $\sigma X$ are co-extensional, it follows that $\epsilon \sigma X = \epsilon X$. Moreover, $\varphi (\sigma X, \sigma \epsilon X, P)$ so $\sigma \epsilon X = \epsilon \sigma X = \epsilon X$---contradiction.\end{proof}
The Axiom of Countable Choice (AC$_\omega$), which states that every countable family of non-empty sets admit a choice function, implies that every infinite set has an countably infinite subset.
\begin{theorem}\label{thm:BLVequivalence1}
Let $U$ be a model of $\ZFUR$ + Collection + AC$_\omega$ equipped with definable classes. The following are equivalent.
\begin{enumerate}
\item $U \models \A$ is a set.
\item $U$ has a definable map that fulfills Basic Law V.
\item $U$ has a parametrically definable map that fulfills Basic Law V.
\end{enumerate}
\end{theorem}
\begin{proof}
(1) $\rightarrow$ (2) follows from Theorem \ref{thm:ASet->BLV}. It remains to prove (3) $\rightarrow$ (1). By Lemma \ref{BLV<->Aset}, it suffices to show that over  $\ZFUR$ + Collection + AC$_\omega$, the principle (*) in Lemma \ref{BLV<->Aset} holds if the urelements do not form a set. So suppose that $\A$ is a proper class and consider any $A \subseteq \A$. Since we have for every $n \in \omega$, there is a $B \subseteq \A$ with $B \cap A = \emptyset$ and $|B| = n$, it follows from Collection that there is an infinite $C \subseteq \A$ that is disjoint from $A$, which has  a countably infinite subset by AC$_\omega$.\end{proof}

The equivalence above holds without Collection if we assume a stronger choice principle. The Axiom of Dependent Choice (DC), which is stronger than AC$_\omega$, states that for every relation on a set without terminal nodes, there is an infinite sequence threading the relation. The \textit{Dependence Choice Scheme }is a class version of the DC.
\begin{itemize}
    \item [] (DC-scheme) If for every $x$ there is some $y$ such that $\varphi(x, y, u)$, then for every $p$ there is an infinite sequence $s$ such that $s(0) = p$ and $\varphi(s(n), s(n+1), u)$ for every $n<\omega$.
\end{itemize}
$\ZFCUR$ + DC-scheme does not prove Collection \cite[Theorem 2.17]{YAO_2024}.

\begin{theorem}\label{thm:BLVequivalence2}
Let $U$ be a model of $\ZFUR$ + DC-scheme equipped with definable classes. The follow are equivalent.
\begin{enumerate}
\item $U \models \A$ is a set.
\item $U$ has a definable map that fulfills Basic Law V.
\item $U$ has a parametrically definable map that fulfills Basic Law V.
\end{enumerate}
\end{theorem}
\begin{proof}
Again, we show that (3) $\rightarrow$ (1) by observing that over  $\ZFUR$ + DC-scheme the principle (*) in Lemma \ref{BLV<->Aset} holds if the urelements do not form a set. Assume that $\A$ is a proper class and let $A \subseteq \A$. Since for every $x$ there is some $y$ such that $ker(x) \subsetneq ker(y)$ and $ker(y) - (A \cup ker(x)) \neq \emptyset$, by the DC-scheme there is an infinite sequence $\<s_n : n< \omega>$ such that $ker(s_n) \subsetneq ker(s_{n+1})$ and $ker(s_{n+1}) - (A \cup ker(s_n)) \neq \emptyset$ for every $n < \omega$. So $ker(s) - A$ is an infinite set of urelements disjoint from $A$, which has a countably infinite subset by the DC-scheme.
\end{proof}
\noindent Note that the argument above also shows that $\ZFCUR$ does not prove the DC-scheme, because (*) fails in the finite-kernel model.

\textit{Plenitude} is the axiom that for every well-ordered cardinal $\kappa$, there is a set of urelements of size $\kappa$. Next, we show that Plenitude is inconsistent with Basic Law V without using any choice principle. We also note that over $\ZFUR$, Plenitude does not imply either Collection or the DC-scheme (\cite[Theorem 17 and Theorem 36]{yao2023set}) although it implies both of them assuming AC (\cite[Theorem 2.5]{YAO_2024}).

\begin{theorem}\label{plenitude->notBLV}
No model of $\ZFUR$ + Plenitude, equipped with definable classes, has a parametrically definable map that fulfills Basic Law V.
\end{theorem}
\begin{proof}
Suppose \textit{for reductio} that in some model $U$ of $\ZFUR$ + Plenitude we can define an extension-assignment map $\epsilon$ using some parameter $P$ defined by some first-order object $z$ in $U$. Note that as in ZF, $\ZFUR$ proves that every set $x$ has a \textit{Hartogs number}, $\aleph(x)$, which is the least ordinal that does not inject into $x$. Let $\kappa = \aleph (ker(z))$ and fix some set of urelements $C$ of size $\kappa$. It follows that $C$ must have a subset $A$ of size $\kappa$ that is disjoint from $ker(z)$. We can then find some $D \subseteq \A$ of size $\kappa^+$, which will have a subset $A'$ of size $\kappa$ that is disjoint from $A \cup ker(z)$. Define $X = \{B \subseteq \A : B - A \text{ is finite}\}$. As in the proof of Claim \ref{claim1.1}, $ker(\epsilon X) \cap A$ cannot be empty since otherwise we can swap $A$ and $A'$ to get a contradiction. But then there is an automorphism $\sigma$ that swaps some urelement $a \in ker(\epsilon X) \cap A$ with some urelement $b \notin ker(\epsilon X)\cup ker(z)$ while fixing $P$. Consequently, $\sigma \epsilon X \neq \epsilon X$, and $\epsilon \sigma X = \epsilon X$ because $X$ and $\sigma X$ are co-extensional. However,  $\sigma$ fixes $P$ so $\sigma \epsilon X = \epsilon \sigma X$---contradiction.\end{proof}

While Hamkins' Theorem \ref{thm:ASet->BLV} suggests that Basic Law V is a natural principle in the context of ZF, the results in this section indicate otherwise: when we quantify over non-sets, Basic Law V entails a specific ontological commitment. For example, Basic Law V implies that there cannot be unboundedly many urelements and that the urelements must form a set when, say, both Collection and AC$_\omega$ hold. 

It is important to note that restricting to \textit{first-order definable classes} does not weaken our argument against Basic Law V. This is because Russell's Paradox shows that Basic Law V is inconsistent if we consider \textit{all} classes. Specifically, the class $R = \{x \mid \exists X (x = \epsilon X \land x \notin X)\}$ cannot have an extension. A related issue is that Basic Law V seems to hold if there exists a \textit{global well-ordering} $<$, since we can then assign every class $\{x \mid \varphi(x)\}$ to the $<$-least ordered-pair $\<\ulcorner \psi_n \urcorner,  u>$ as in the proof of Theorem \ref{thm:ASet->BLV}. The status of the Axiom of Choice and its second-order variants in urelement set theory will be discussed later. However, the problem with this argument is that there is no \textit{first-order definable} global well-ordering if the urelements form a proper class. For, if $\varphi(x ,y ,p)$ defines a global well-ordering $<$ and the urelements do not form a set, then swapping two urelements $a$ and $b$ outside $ker(p)$ such that $a < b$ would yield $b < a$, contradicting the assumption of a well-ordering.

\section{Hume's Principle in Urelement Set Theory}\label{section:HPforSets}
A map $x \mapsto \# x$ from sets to first-order objects fulfills \textit{Hume's Principle} if for any sets $x$ and $y$, $x$ is equinumerous to $y$ ($x \sim y$) just in case $\# x = \# y$. We say that Hume's Principle holds for sets when such a map exists. In this section, we investigate the independence of Hume's Principle for sets in ZF with urelements. The philosophical implications of this independence will be discussed at the end.

To begin with, $\ZFUR$ is consistent with Hume's Principle for sets. When the urelements form a set, we can again use Scott's trick to map each set $x$ to $\# x = \{y : y \sim x \land \forall z (z \sim x \rightarrow \rho (z) \leq \rho(y))\}$.
And no matter if the urelements form a set, AC implies that Hume's Principle holds for sets since we can map each set $x$ to the least ordinal to which $x$ is equinumerous.

However, we note that neither AC nor urelements' forming a set is necessary for Hume's Principle to hold for sets. A set $A$ of urelements is said to be \textit{universal} if every set $x$ is equinumerous with some set $y$ in $V(A)$.

\begin{prop}
If there is a universal set of urelements, then Hume's Principle holds for sets.
\end{prop}
\begin{proof}
We simply use Scott's trick within in $V(A)$. For every set $x$, let $$\# x = \{y \in V(A) : y \sim x \land \forall z \in V(A) (z \sim x \rightarrow \rho (z) \leq \rho(y))\}.$$
$\# x$ is a well defined set because $A$ is universal. It is then routine to check that the map $x \mapsto \# x$ fulfills Hume's Principle.
\end{proof}
\noindent There are models of $\ZFUR$ + Collection + $\neg$AC in which the urelements form a proper class but a universal set of urelements exists (see, e.g., the model $W$ in \cite[Theorem 36]{yao2023set}). Although it is not known whether Hume's Principle for sets implies the existence of a universal set of urelements, we note that the following weaker implication holds, the proof of which will be useful later.
\begin{prop}\label{prop:HP->AlmostUniverseSet}
If Hume's Principle holds for sets, then there is a set of urelements $A$ such that $\# x \in V(A)$ for every set $x$.\footnote{Thus, with an additional requirement that $x \sim \# x$ for every set $x$ (i.e., every set is equinumerous with its cardinality), Hume's Principle implies that a universal set of urelements exists.It is proved in \cite{pincus1974cardinal} that this stronger form of Hume's Principle is not provable in ZF and does not imply AC.}
\end{prop}
\begin{proof}
We may assume the urelements form a proper class. Suppose that some formula $\varphi(x, \# x, u)$ defines a map that fulfills Hume's Principle using a parameter $u$. We show that $\# x \in V(ker(u))$ for every set $x$. Suppose not. Fix some $a \in ker(\# x) - ker(u)$ and some $b \notin ker(u) \cup ker(\# x)$. Let $\pi$ be the automorphism that only swaps $a$ and $b$. So $\pi \# x \neq \# x$. Then we have $\varphi(\pi x, \pi \# x, u)$ so $\pi \# x = \# \pi x$. But $\pi x \sim x$, so $\# \pi x = \# x$---contradiction. \end{proof}

We now move on to the independence of Hume's Principle for sets in $\ZFUR$. A known result is the following.
\begin{theorem}[Gauntt \cite{gauntt1967undefinability}; L\'evy \cite{levy1969definability}]
 There are models of $\ZFUR$ in which Hume's Principle does not hold for sets. \qed
\end{theorem}
\noindent As we noted in Section \ref{section:AxiomsinZFU}, since $\ZFUR$ is a rather weak theory the independence of Hume's Principle will be more robust if the counter-models of Hume's Principle also satisfy stronger axioms such as RP. In other words, given the theorem by Gauntt and L\'evy, it is natural to ask whether Hume's Principle depends upon the axiomatization of urelement set theory (in the choiceless context, of course). To address this question, we shall construct a model where Hume's Principle fails for sets but RP holds. Thus, Hume's Principle for sets is independent of ZF with urelements regardless of one's axiomatization of the theory.

Our proof that RP holds in the model will again appeal to two aforementioned axioms, Collection and the DC-scheme, due to the following observation made in \cite[page 397]{Gitman2016-GITWIT} (though in a different context).
\begin{theorem}\label{collection+dc->rp}
$\ZFUR \vdash $ Collection $\land$ DC-scheme $\rightarrow$ RP. \qed
\end{theorem}
\noindent In fact, our result will also clarify the relationship between Hume's Principle and AC by showing that \textit{no fragment} of AC suffices for Hume's Principle. By ``fragment of AC'', we mean the DC$_\kappa$-hierarchy. We shall review some facts about the DC$_\kappa$-hierarchy that will be useful. To begin with, for every infinite well-ordered cardinal $\kappa$, DC$_\kappa$ is the following axiom generalizing DC.
\begin{itemize}
    \item [] (DC$_\kappa$) For every set $x$ and relation $r\subseteq x^{<\kappa} \times x$, if for every $s \in x^{<\kappa}$, there is some $w \in x$ such that $\<s, w> \in r$, then there is an $f: \kappa \rightarrow x$ such that $\<f\restriction\alpha, f(\alpha)> \in r$ for all $\alpha < \kappa$.
\end{itemize}
As in ZF, $\ZFUR$ proves that AC is equivalent to $\forall \kappa$ DC$_\kappa$, so each DC$_\kappa$ is a fragment AC. Moreover, for each $\kappa$, we can formulate the DC$_\kappa$-scheme as the class version of DC$_\kappa$.
\begin{itemize}
        \item[]  (DC$_\kappa$-scheme) If for every $x$ there is some $y$ such that $\varphi(x, y, u)$, then there is some function $f : \kappa\rightarrow U$ such that $\varphi(f\restriction \alpha, f(\alpha), u)$ for every $\alpha <\kappa$.
\end{itemize}
By a standard argument $\ZFCUR$ + $\A$ is a set proves that the DC$_\kappa$-scheme for every $\kappa$ (\cite[Lemma 23]{yao2023set}). If $\lambda < \kappa$, the DC$_\kappa$-scheme implies the DC$_\lambda$-scheme; and the DC$_\omega$-scheme is a reformulation of the DC-scheme. Moreover, the DC$_\kappa$-scheme is equivalent to the following (\cite[Proposition 13]{yao2023set}).
\begin{itemize}
\item [] For every definable class $X$, if for every $s \in X^{<\kappa}$ there is some $y \in X$ with $\varphi(x, y, u)$, then there is some function $f : \kappa \rightarrow X$ such that $\varphi(f\restriction \alpha, f(\alpha), u)$ for every $\alpha < \kappa$.
\end{itemize} 

\begin{theorem}\label{thm:HPindependence}
Let $\kappa$ be any well-ordered infinite cardinal. There is a model of $\ZFUR$ + RP + DC$_\kappa$-scheme in which Hume's Principle does not hold for sets. In fact, in this model no parametrically definable map fulfills Hume's Principle.
\end{theorem}

\begin{proof}
The first step is to build a \textit{permutation model} that violates AC. We assume some familiarity with permutation models and only provide some necessary details of the construction. For a detailed presentation of this topic, see \cite{jech2008axiom}. Fix some infinite cardinal $\kappa$. We start with a model $U$ of $\ZFCUR$ in which $\A$ is a set of size $\kappa^+$. Enumerate $\A$ with $\kappa^+ \times \kappa^+$, i.e., $\A = \bigcup_{\alpha < \kappa^+} A_\alpha$, where each row $A_\alpha$ has size $\kappa^+$. Let $\G$ be the group of permutations of $\A$ such that for every $\pi \in \G$ and $\alpha < \kappa^+$, $\pi A_\alpha = A_\alpha$. Let $I = \{E \subseteq \A : |E \cap A_\alpha| < \kappa^+ \text{ for each } \alpha < \kappa^+\}$ be an ideal of $\A$. For every $x$, define $sym(x) = \{\pi \in \G_A : \pi x = x\}$; if $x$ is a set, define $fix(x) = \{\pi \in \G_A : \pi y = y \text{ for all } y\in x \}$. We say an object is \textit{symmetric} if there some $E \in I$, called a \textit{support} of $x$, such that $fix(E) \subseteq sym(x)$. The permutation model $W$ is the class of all \textit{hereditarily symmetric} objects, i.e., $W = \{x \in U : x \text{ is symmetric} \land x \subseteq W\}$, which is a model of $\ZFUR$ (see \cite[Theorem 4.1]{jech2008axiom}, or \cite[Theorem 33]{yao2023set} for a more general proof).

Note that $\A \in W$ so there is only a set of urelements in $W$ and consequently Hume's Principle still holds in $W$. The next step is to produce a model with proper-class many urelements in which Hume's Principle fails. In $W$, define a class of sets $\I$ of urelements such that
 $$\I = \{B \subseteq \A : B \text{ is a subset of } \kappa\text{-many } A_\alpha\}.$$
 \noindent Let $W^\I = \{ x \in W : ker(x) \in \I\}$. $\I$ is a class \textit{ideal} of $\A$ which picks out some ``small'' sets of urelements, while $W^\I$ is the class of all objects in $W$ whose kernel is small in the sense of $\I$. By \cite[Theorem 26]{yao2023set}, $W^\I \models \ZFUR$ and there is a proper class of urelements in $W^\I$. 
\begin{lemma}\label{lemma:WImodelsDCK}
$W^\I \models $ DC$_\kappa$-scheme.
\end{lemma}
\begin{proof}
First observe that $W^\I$ is closed-under $\kappa$-sequences. Consider any $s : \kappa \rightarrow W^\I$. In $U$, for each $\alpha < \kappa$ choose some $E_\alpha \in I$ to be a support of $s_\alpha$. Then $\bigcup_{\alpha< \kappa}E_\alpha \in I$, which is a support of $s$. And since $ker(s) = \bigcup_{\alpha < \kappa} ker(s_\alpha)$, it follows that $ker(s)$ is also contained in a $\kappa$-block of $A_\alpha$. Thus, $s \in W^\I$. Now suppose that $W^\I \models \forall x \exists y \varphi(x, y, u)$, where $u \in W^\I$. Then for every $s \in (W^\I)^{<\kappa}$ there is some $y \in W^\I$ such that $\varphi^{W^\I}(s, y, u)$. By the DC$_\kappa$-scheme in $U$, there is an $f: \kappa \rightarrow W^\I$ with $\varphi^{W^\I}(f\restriction\alpha, f_\alpha, u)$ for every $\alpha < \kappa$. Since $f \in W^\I$, it follows that $W^\I \models $ DC$_\kappa$-scheme.
\end{proof}

Next we show that RP holds in $W^\I$. By Theorem \ref{collection+dc->rp} and Lemma \ref{lemma:WImodelsDCK}, to show $W^\I \models$ RP it suffices to check that $W^\I \models$ Collection. A permutation $\sigma$ of $\A$ in $U$ is said to be \textit{row-swapping} if for every $\alpha < \kappa^+$, $\sigma A_{\alpha} = A_{\beta}$ for some $\beta < \kappa^+$ (consequently,  for every $\alpha$, $A_{\alpha} = \sigma A_{\beta}$ from some $\beta$).
 \begin{lemma}\label{Wfix}
Every row-swappoing $\sigma$ in $U$ is an automorphism of $W^\I$. 
\end{lemma}
\begin{proof}
Let $\sigma$ be a row-swappoing permutation of $\A$. It suffices to show that $\mathcal{G}_\A, \ I$, and $\I$ are all fixed by $\sigma$. If $\pi \in \G_\A$, then for every $A_\alpha$, since $A_\alpha = \sigma A_\beta$ for some $\beta$ and $\pi A_\beta = A_\beta$, by automorphism it follows that $(\sigma \pi) (\sigma A_\beta) = \sigma A_\beta$; so $(\sigma \pi )A_\alpha = A_\alpha$ and hence $\sigma \pi \in \mathcal{G}_\A$ (note here $\sigma \pi$ is not $\sigma \circ \pi$ but $\{\langle \sigma a, \sigma (\pi a) \rangle : a \in \A \rangle\})$. This shows that $\sigma \mathcal{G}_\A = \mathcal{G}_\A$. If $E \in I$, then for every $A_\alpha$, since $A_\alpha = \sigma A_\beta$ for some $\beta$ and $E \cap A_\beta$ has size $< \kappa^+$, $\sigma E \cap A_\alpha$ has size $< \kappa^+$ and hence $\sigma E \in I$. Therefore, $\sigma I= I$. Similarly, if $B$ is contained in $\kappa$-many $A_\alpha$, then so is $\sigma B $. Hence, $\sigma \I= \I$ and the lemma is proved.
\end{proof}

\begin{lemma}\label{lemma:WImodelsRP}
$W^\I \models$ RP.
\end{lemma}
\begin{proof}
We show that $W^\I \models$ Collection. Suppose that $W^\I\models \forall x \in w \exists y \varphi(x, y, u)$ for some $w, u \in W^\I$. Let $A=\bigcup_{\alpha < \kappa} A_{\alpha}$ be a $\kappa$-block containing $ker(w) \cup ker(u)$ and $B = \bigcup_{\alpha < \kappa} B_{\alpha}$ be another $\kappa$-block that is disjoint from $A$. It is enough to show that $W^\I\models \forall x \in w \exists y \in V(A \cup B) \varphi(x, y, u)$ because then a sufficiently tall $V_\alpha (A \cup B)$ will be a desired collection set. 

Consider any $x \in w$ and $y \in W^\I$ such that $W^\I \models \varphi(x, y, u)$. Let $C = \bigcup_{\alpha < \kappa} C_{\alpha}$ be another $\kappa$-block containing $ker(y)-A$. 
\begin{claim}
In $U$ there is a row-swapping $\sigma$ such that $\sigma(B \cup C) = B$.
\end{claim}
\begin{claimproof}
Fix a $\kappa$-block $D = \bigcup_{\alpha < \kappa} D_{\alpha}$ disjoint from $A \cup B \cup C$. Split $B$ into two pair-wise disjoint $\kappa$-blocks $B^1= \bigcup_{\alpha < \kappa} B^1_{\alpha}$ and $B^2= \bigcup_{\alpha < \kappa} B^2_{\alpha}$, where each $B^1_{\alpha}$ (and $B^2_{\alpha}$) is some row $A_\alpha$. Then split $D$ into two pair-wise disjoint $\kappa$-blocks $D^1= \bigcup_{\alpha < \kappa} D^1_{\alpha}$ and $D^2= \bigcup_{\alpha < \kappa} D^2_{\alpha}$ in the same way. In $U$ we can define a row-swapping permutation $\sigma$ of $\A$ as follows. For each $\alpha < \kappa^+$, let $\sigma B_\alpha = B^1_\alpha$, $\sigma C_\alpha = B^2_\alpha$, $\sigma D^1_\alpha = C_\alpha$ and $\sigma D^2_\alpha = D_\alpha$. It is clear that $\sigma(B \cup C) = B$.\end{claimproof}
\vspace{4pt}

\noindent By Lemma \ref{Wfix} $\sigma$ is an automorphism of $W^\I$ fixing $x$ and $u$, so $W^\I \models \varphi(x, \sigma y, u) \land \sigma y \in  V(A \cup B)$, which proves the lemma.\end{proof}

Finally, we turn to the failure of Hume's Principle in $W^\I$.
Suppose \textit{for reductio} that $\varphi(x, \# x, u)$ defines a cardinality-assignment map with some parameter $u \in W^\I$. It follows that there must be two $A_\alpha$ and $A_\beta$ that are disjoint from $ker(u)$. 

\begin{claim}\label{alphabetanotequi}
    $W^\I \models A_\alpha \nsim A_\beta$.
\end{claim}
\begin{claimproof}
Suppose \textit{for reductio} that $f$ is an injection from $A_\alpha$ to $A_\beta$ in $W^\I$. Let $E \in I$ be a support of $f$. Then there are two urelements $a, b \in A_\alpha - E$. Let $\pi \in fix(E)$ swap only $a$ and $b$. Then $f(b) = \pi f(a) = f(a)$, contradicting the assumption that $f$ is an injection. \end{claimproof}
\vspace{4pt}

\noindent By Lemma \ref{Wfix}, in $U$ there is an automorphism $\sigma$ of $W^\I$ that swaps only $A_\alpha$ and $A_\beta$. Since $\sigma$ point-wise fixes $ker(u)$ and $ker(\# A_\alpha) \subseteq ker(u)$ by the proof of Proposition \ref{prop:HP->AlmostUniverseSet}, $\sigma$ point-wise fixes $ker(\# A_\alpha)$. It follows that $W^\I \models \varphi(A_\beta, \# A_\alpha, u)$ so $\# A_\beta = \# A_\alpha$. Hence, $W^\I \models A_\alpha \sim A_\beta$, contradicting Claim \ref{alphabetanotequi}. This completes the proof of the theorem.
\end{proof}
We now consider the philosophical implications of Theorem \ref{thm:HPindependence} for Hume's Principle. Since we are working in ZF with urelements, we shall assume a \textit{predicativist} view of classes—that is, all classes are first-order definable. Hume's Principle in impredicative class theories will be discussed in Section \ref{section:2ndOrderHP}. As noted in the introduction, one desideratum for Hume's Principle to be regarded as a fundamental principle of cardinality is its \textit{universality}, meaning that the opening quantifiers in Hume's Principle must range over concepts of any kind. Given this requirement, models of urelement set theory with definable classes are precisely the context in which Hume's Principle should be tested. Consequently, the model $W^\I$ in Theorem \ref{thm:HPindependence} provides evidence against the universality of Hume's Principle, assuming it is regarded as a natural model of set theory with urelements.  

One possible objection to the naturalness of $W^\I$ comes from proponents of the Axiom of Choice, who may argue that $W^\I$ is unnatural because mathematics without AC can exhibit various pathologies. While it is well known that AC itself leads to certain paradoxical consequences, let us grant, for the sake of argument, that AC is an overall desirable axiom in pure set theory. In response, we note that AC does hold in the class of \textit{pure sets} ($V$) of $W^\I$, as no pure sets were eliminated in the process of constructing the model. Therefore, since all mathematical objects can be reconstructed within the universe of pure sets, no mathematical pathologies that would typically arise in the absence of AC occur in $W^\I$. Moreover, as we observed earlier, AC asserts more than what is necessary for Hume's Principle to hold: the existence of a universal set of urelements suffices. Thus, there seems to be no compelling reason to consider the claim ``\textit{every set} is well-orderable'' a natural axiom in urelement set theory.

\section{Hume's Principle in Urelement Class Theory}\label{section:2ndOrderHP}
A map $X \mapsto \# X$ from classes to first-order objects fulfills Hume's Principle if for any classes $X$ and $Y$, $X \sim Y$ just in case $\# X = \# Y$.  When a map as such exists, we say that Hume's Principle holds \textit{for classes}. Unlike Basic Law V, Hume's Principle for classes is consistent with impredicative class theories such as the Kelley-Morse class theory (KM). In particular, Hume's Principle will hold for classes if von Neumann's \textit{Limitation of Size} holds, which states that all proper classes are equinumerous. This is because Limitation of Size implies the Axiom of Global Well-Ordering, which asserts the existence of a well-ordering of the universe. So with Limitation of Size, we can define, without using parameters, the following map $X \mapsto \# X$ that fulfills Hume's Principle.
\begin{equation*}
    \# X = 
    \begin{cases*}
       \{\{0\}\}    & if $X \text{ is a proper class}$ \\
       \alpha & $\alpha$ is the least ordinal such that $X \sim \alpha$
    \end{cases*}
\end{equation*}
And unlike Hume's Principle for sets, without assuming AC Hume's Principle can fail for classes even when there are no urelements (see \cite[Theorem 2]{Roberts2023-ROBHPB}).

Let KMCU\footnote{We do not distinguish different axiomatizations here because with Global Well-Ordering, Replacement, Collection and RP are all equivalent. Yet, without Global Well-Ordering the situation is similar to ZFC with urelements. See \cite[Chapter 4]{yao2023set}.} be KMU plus Global Well-Ordering. In KMCU, let $W$ be a global well-ordering on $U$. Then we can use $W$ as a class parameter to define a map $X \mapsto \# X$ as follows.\footnote{We are grateful to an anonymous referee for pointing this out.}
\begin{equation*}
    \# X = 
    \begin{cases*}
       0    & if $X \sim U$ \\
       \<0, y> & $X \not\sim U$ and $y$ is the $W$-least object such that $X \sim W_y$ 
    \end{cases*}
\end{equation*}
where $W_y = \{x \mid x W y\}$ is the initial segment of $W$ below $y$. This is a well-defined map due to the fact that for every class $X$, the well-order $W\restriction X$ is either isomorphic to $W$ or isomorphic to $W_x$ for some $x$.\footnote{A sketch of proof of this fact. Define class $$F=\{\<w, y> \in X \times U \mid (W\restriction X)_w \text{ is isomorphic to } W_y \}.$$ $F$ is an isomorphism between $dom(F)$ and $ran(F)$. As both $dom(F)$ and $ran(F)$ are closed downward with respect to $W$, either $dom(F) = X$ or $ran(F) = U$. Moreover, if $ran(F) = U$, then $dom(F) = X$ since otherwise $W$ will embed into a bounded proper segment of itself, which is impossible.} Thus, if $X \not\sim U$, it must be isomorphic and hence equinumerous to some $W_x$. It is then routine to check that $X \mapsto \# X$ so defined fullfills Hume's Principle.

A natural question is whether such a map $X \mapsto \# X$ is definable \textit{without parameters} in KMCU. A standard argument shows that Global Well-Ordering implies Limitation of Size with the assumption that the urelements form a set (\cite[Proposition 98]{yao2023set}), so the answer to the question is positive if the urelements are ``few'' compared to sets. It is shown in \cite{howard1978independence} that KMCU by itself does not prove Limitation of Size. The next theorem shows that Global Well-Ordering is not sufficient for a parameter-free definable map that fulfills Hume's Principle when the urelements are ``many'' compared to the pure sets.

\begin{theorem}\label{KMCUnvdashHP}
Assume the consistency of ZFC $+$ there is an inaccessible cardinal. There is a model of KMCU which has no definable map that fulfills Hume's Principle for classes.
\end{theorem}
\begin{proof}
Let $V$ be a model of ZFC $+$ an inaccessible cardinal $\kappa$. We first construct a model $U$ of $\ZFCUR$ + Plenitude by treating copies of ordinals as urelements. In particular, we define $V \llbracket Ord \rrbracket$ in $V$ by recursion as follows.
\begin{itemize}
    \item [] $ V \llbracket Ord \rrbracket = (\{0\} \times Ord) \cup \{\barx \in V : \exists x (\barx = \<1 , x> \land x \subseteq V \llbracket Ord \rrbracket)\}.$
\end{itemize}
For every $\barx, \bary \in V \llbracket Ord \rrbracket$,
\begin{itemize}
    \item [] $\barx \ \bar{\in} \ \bary \text{ if and only if } \exists y (\bary = \<1, y> \land  \barx  \in y);$
    \item [] $\bar{\A}(\barx) \text{ if and only if } \barx \in \{0\} \times Ord.$
\end{itemize}
Let $U$ denote the model $\<V \llbracket Ord \rrbracket,\ \bar{\A},\ \bar{\in}>$ for the language of urelement set theory. By \cite[Theorem 7 and 9]{yao2023set}, $U \models \ZFCUR$ + Plenitude. Moreover, the class of pure sets in $U$ is isomorphic to $V$ (\cite[Lemma 8]{yao2023set}) so $U$ also contains (a copy of) $\kappa$ as an inaccessible cardinal.

In  $U$, let $\lambda = \aleph_{\kappa^+}$ and $A$ be a set of urelements of size $\lambda$. Define $$U_\kappa(A) = \bigcup_{B \in P_\kappa (A)} V_\kappa (B),$$ where $P_\kappa (A)$ is the set of subsets of $A$ with size $< \kappa$. By \cite[Lemma 4.2]{yao2022reflection}, $U_\kappa(A)$, equipped with all of its subsets, is a model of KMCU. Now suppose \textit{for reductio} that $\varphi(X, \# X)$ defines a cardinality-assignment function in $U_\kappa (A)$. 
\begin{claim}
For every class $X$ of $U_\kappa (A)$, $\# X$ is a pure set.
\end{claim}
\begin{claimproof}
Suppose $\# X$ is not pure for some class $X$ of $U_\kappa (A)$. Then in $U_\kappa(A)$ there will be an automorphism $\pi$ swapping some $a \in ker(\#X)$ with some $b \notin ker(\#X)$, and so $\pi \# X \neq \# X$. Since $\pi X \sim X$, it follows that $\# X = \# \pi X = \pi \# X$, yielding a contradiction.\end{claimproof}
\vspace{6pt}

\noindent Since there are at least $\kappa^+$-many cardinalities for the proper classes of $U_\kappa(A)$, $\varphi$ thus defines an injection from $\kappa^+$ into $V^{U_\kappa(A)}$, namley, $V_\kappa$. This is impossible as $V_\kappa$ has size $\kappa$.\end{proof}

\noindent In fact, this argument shows that in $U_\kappa(A)$ there is no definable cardinality-assignment map which only uses first-order objects in $U_\kappa(A)$ as parameters. For, if  $\varphi(X, \# X, u)$ defines a cardinality-assignment map, where $u \in U_\kappa(A)$, then the same argument would show that $ker(\# X) \subseteq \ker (u)$ for every class $X$. This means that we could inject $\kappa^+$-many cardinalities into $V_\kappa (ker(u))$, but $V_\kappa (ker(u))$ only has size $\kappa$ as $|ker(u)| < \kappa$---contradiction.

\section*{Acknowledgement}
The author would like to thank the audience at the Logic Colloquium of Wuhan University and two anonymous referees for their helpful and detailed comments, with special thanks to one referee who points out that there is always a parametrically definable map in KMCU that fulfills Hume's Principle.

\section*{Funding}
The author was supported by NSFC No. 12401001 and the Fundamental Research Funds for the Central Universities, Peking University.
\printbibliography

@article{tait2005constructing,
  title={Constructing cardinals from below},
  author={Tait, William W},
  journal={The Provenance of Pure Reason: essays in the philosophy of mathematics and its history},
  pages={133--154},
  year={2005}
}

@article{howard1978independence,
  title={Independence results for class forms of the axiom of choice},
  author={Howard, Paul E and Rubin, Arthur L and Rubin, Jean E},
  journal={Journal of Symbolic Logic},
  pages={673--684},
  year={1978},
  publisher={JSTOR}
}

@incollection{bernays1976problem,
  title={On the problem of schemata of infinity in axiomatic set theory},
  author={Bernays, Paul},
  booktitle={Studies in Logic and the Foundations of Mathematics},
  volume={84},
  pages={121--172},
  year={1976},
  publisher={Elsevier}
}

@book{jech2008axiom,
  title={The axiom of choice},
  author={Jech, Thomas J},
  year={2008},
  publisher={Courier Corporation}
}

@article{Gitman2016-GITWIT,
	journal = {Mathematical Logic Quarterly},
	volume = {62},
	doi = {10.1002/malq.201500019},
	pages = {391--406},
	author = {Victoria Gitman and Joel David Hamkins and Thomas A. Johnstone},
	title = {What is the Theory ZFC Without Power Set?},
	number = {4-5},
	year = {2016}
}

@article{yao2022reflection,
  title={Reflection principles and second-order choice principles with urelements},
  author={Yao, Bokai},
  journal={Annals of Pure and Applied Logic},
  volume={173},
  number={4},
  pages={103073},
  year={2022},
  publisher={Elsevier}
}

@incollection{levy1969definability,
  title={The definability of cardinal numbers},
  author={L{\'e}vy, Azriel},
  booktitle={Foundations of Mathematics},
  pages={15--38},
  year={1969},
  publisher={Springer}
}

@article{levy1961principles,
  title={Principles of partial reflection in the set theories of Zermelo and Ackermann.},
  author={L{\'e}vy, Azriel and Vaught, Robert},
  journal={Pacific journal of mathematics},
  volume={11},
  number={3},
  pages={1045--1062},
  year={1961},
  publisher={Pacific Journal of Mathematics, A Non-profit Corporation}
}

@book{yao2023set,
  title={Set theory with urelements},
  author={Yao, Bokai},
  year={2023},
  eprint={2303.14274},
    archivePrefix={arXiv},
  publisher={University of Notre Dame}
}

@misc{hamkins2022fregean,
      title={Fregean abstraction in Zermelo-Fraenkel set theory: a deflationary account}, 
      author={Joel David Hamkins},
      year={2022},
      eprint={2209.07845},
      archivePrefix={arXiv},
      primaryClass={math.LO}
}

@article{Bell1994-BELFEO,
	author = {John L. Bell},
	doi = {10.1002/malq.19940400104},
	journal = {Mathematical Logic Quarterly},
	number = {1},
	pages = {27--30},
	publisher = {Wiley-Blackwell},
	title = {Fregean Extensions of First-Order Theories},
	volume = {40},
	year = {1994}
}

@article{Burgess1998-JOHOAC,
	author = {John P. Burgess},
	doi = {10.1305/ndjfl/1039293068},
	journal = {Notre Dame Journal of Formal Logic},
	number = {2},
	pages = {274--278},
	publisher = {Duke University Press},
	title = {On a Consistent Subsystem of Frege's Grundgesetze},
	volume = {39},
	year = {1998}
}

@article{Parsons1987-TEROTC,
	author = {Terence Parsons},
	doi = {10.1305/ndjfl/1093636853},
	journal = {Notre Dame Journal of Formal Logic},
	number = {1},
	pages = {161--168},
	publisher = {Duke University Press},
	title = {On the Consistency of the First-Order Portion of Frege's Logical System},
	volume = {28},
	year = {1987}
}

@article{gauntt1967undefinability,
  title={Undefinability of cardinality},
  author={Gauntt, Robert J},
  journal={Lectures notes prepared in connection with the Summer Institute on Axiomatic Set Theory held at University of California, Los Angeles, IV-M},
  year={1967}
}

@book{Fine2002-FINTLO-4,
	address = {New York},
	author = {Kit Fine},
	editor = {Matthias Schirn},
	publisher = {Oxford University Press},
	title = {The Limits of Abstraction},
	year = {2002}
}

@article{Shapiro2005-SHASAA-12,
	author = {Stewart Shapiro},
	doi = {10.1007/s11098-004-4923-9},
	journal = {Philosophical Studies},
	number = {3},
	pages = {315--332},
	publisher = {Springer},
	title = {Sets and Abstracts ? Discussion},
	volume = {122},
	year = {2005}
}

@article{Roberts2023-ROBHPB,
	author = {Sam Roberts and Stewart Shapiro},
	doi = {10.1017/s1755020322000144},
	journal = {Review of Symbolic Logic},
	number = {4},
	pages = {1158--1176},
	title = {Hume?s Principle, Bad Company, and the Axiom of Choice},
	volume = {16},
	year = {2023}
}

@article{Boolos1989-BOOIA,
	author = {George Boolos},
	doi = {10.5840/philtopics19891721},
	journal = {Philosophical Topics},
	number = {2},
	pages = {5--21},
	publisher = {University of Arkansas Press},
	title = {Iteration Again},
	volume = {17},
	year = {1989}
}

@article{Shapiro1999-RICNVZ,
	author = {Stewart Shapiro and Alan Weir},
	doi = {10.1093/philmat/7.3.293},
	journal = {Philosophia Mathematica},
	number = {3},
	pages = {293--321},
	publisher = {Oxford University Press},
	title = {New V, Zf and Abstraction},
	volume = {7},
	year = {1999}
}

@article{Shapiro2008-SHAFMZ-2,
	author = {Stewart Shapiro},
	doi = {10.1017/s1755020308080192},
	journal = {Review of Symbolic Logic},
	number = {2},
	pages = {241--266},
	publisher = {Cambridge University Press},
	title = {Frege Meets Zermelo: A Perspective on Ineffability and Reflection},
	volume = {1},
	year = {2008}
}

@article{YAO_2024, title={AXIOMATIZATION AND FORCING IN SET THEORY WITH URELEMENTS}, DOI={10.1017/jsl.2024.58}, journal={The Journal of Symbolic Logic}, author={YAO, BOKAI}, year={2024}, pages={1–27}}

@book{chang1990model,
  title={Model theory},
  author={Chang, Chen Chung and Keisler, H Jerome},
  volume={73},
  year={1990},
  publisher={Elsevier}
}

@article{bell1977boolean,
  title={Boolean-valued models and independence proofs in set theory},
  author={Bell, John Lane},
  year={1977}
}

@article{pincus1974cardinal,
  title={Cardinal representatives},
  author={Pincus, David},
  journal={Israel Journal of Mathematics},
  volume={18},
  pages={321--344},
  year={1974},
  publisher={Springer}
}

\end{document}